\documentstyle[12pt]{article}                  
\input{amssymb.sty}
\makeatletter
\def\nothing#1{}
\newdimen\earraycolsep
\renewcommand{\thetable}{\arabic{table}}
\renewcommand{\thefigure}{\arabic{figure}}
\renewcommand{\title}[1]{%
  \vspace*{50\p@}%
   {\parindent \z@ \raggedright \reset@font
     \bfseries #1\par
     \nobreak
     \vskip 36\p@
   }}
\def\author#1{{\pretolerance=10000 \raggedright \advance \leftskip by 1in
\noindent #1 \vskip 1pc}}
\def\affiliation#1{{\advance\leftskip by 1in \noindent #1 \vskip -1pc}}
\def\refnote#1{{$^{\hbox{\scriptsize #1}}$}}
\def\affnote#1{\llap{$^{\hbox{\scriptsize #1}}$}}

\renewcommand\section{\@startsection{section}{1}{\z@}{2pc \@plus
       1ex minus .2ex}{1pc \@plus .2ex}{\reset@font
       \normalsize\bfseries\noindent
       {\addtocounter{section}{1}}\Roman{section}\
       {\setcounter{subsection}{0}
       \setcounter{subsubsection}{0}\setcounter{equation}{0}} }}
\renewcommand\subsection{\@startsection{subsection}{2}{\z@}{1pc \@plus 1ex
     minus.2ex}{1pc \@plus .2ex}
     {\reset@font\normalsize\bfseries
     \noindent{\addtocounter{subsection}{1}}%
     {\setcounter{subsubsection}{0}}\Roman{section}.\Roman{subsection}\ }}
\renewcommand\subsubsection{\@startsection{subsubsection}{3}{\parindent}
         {1pc \@plus 1ex minus.2ex}{-0.5em}{\reset@font\normalsize\bfseries%
         {\addtocounter{subsubsection}{1}} \hspace*{.6cm}
         \Roman{section}.\Roman{subsection}.\Roman{subsubsection}
         \hspace*{-7mm}}}

\def\AmS{{\protect\the\textfont2%
         A\kern-.1667em\lower.5ex\hbox{M}\kern-.125emS}}

\def\p@LaTeX{{\family{times}\series{m}\shape{n}\selectfont
L\kern-.36em\raise.3ex\hbox{\scriptsize A}\kern-.15em
T\kern-.1667em\lower.7ex\hbox{E}\kern-.125emX}}

\newlength{\colwidth}

\def\@oddhead{\hfil}
\def\@evenhead{\hfil}
\def\@oddfoot{{\bfseries\hfil\thepage}}
\def\@evenfoot{{\bfseries\thepage\hfil}}
\def\fnum@figure{\footnotesize\raggedright{\bfseries \figurename~\thefigure.}}
\def\fnum@table{\normalsize\raggedright{\bfseries \tablename~\thetable.}}
\long\def\@makecaption#1#2{\vskip 10\p@ {#1 #2\par}}
\long\def\@makefntext#1{\setbox0=\hbox{$\m@th^{\@thefnmark}$}\noindent
\hangindent=\wd0 \box0 #1}
\newbox\@atbox
\long\def\atable#1#2#3{\begin{table}[tbp]\centering\footnotesize
\setbox\@atbox\hbox{#2}
\parbox{\wd\@atbox}{\caption{#1}}\par\smallskip
#2
\par\smallskip\parbox{\wd\@atbox}{\raggedright #3}
\end{table}}

\newtheorem{theorem}{Theorem}

\newtheorem{lemma}[theorem]{Lemma}

\def\Cb{{\mathbb C}}
\def\Hb{{\mathbb H}}
\def\Nb{{\mathbb N}}
\def\Rb{{\mathbb R}}
\def\Tb{{\mathbb T}}
\def\Zb{{\mathbb Z}}

\def\Ac{{\cal A}}
\def\Bc{{\cal B}}
\def\Ac{{\cal A}}
\def\Hc{{\cal H}}
\def\Lc{{\cal L}}

\def\Hc{{\cal H}}

\def\Lc{{\cal L}}

\newcommand{\wt}{\widetilde}

\def\a{\alpha}
\def\b{\beta}
\def\d{\delta}
\def\lb{\lambda}
\def\g{\gamma}
\def\om{\omega}
\def\s{\sigma}
\def\t{\theta}

\def\vp{\varphi}

\def\G{\Gamma}

\def\fl{\forall}
\def\ify{\infty}
\def\lgl{\langle}

\def\ot{\otimes}
\def\ov{\overline}
\def\part{\partial}
\def\rgl{\rangle}
\def\sbs{\subset}

\def\ts{\times}

\def\ra{\rightarrow}
\def\longra{\longrightarrow}

\def\text{\hbox}

\def\Sup{\mathop{\rm Sup}\nolimits}

\def\build#1_#2^#3{\mathrel{
\mathop{\kern 0pt#1}\limits_{#2}^{#3}}}

\def\@nbibitem#1{\noindent \hangindent=2pc \hangafter=1
\refstepcounter{enumi}\hbox to 2pc{\arabic{enumi}.\hfil}%
\immediate\write\@auxout{\string\bibcite{#1}{\arabic{enumi}}}}
\def\numbibliography{%
\section*{REFERENCES}%
\bgroup\footnotesize
\setcounter{enumi}{0}%
\def\newblock{\hskip .11em plus.33em minus.07em}%
\let\bibitem\@nbibitem}
\def\endnumbibliography{\par\egroup}
\makeatother
\def\ra{\rightarrow}
\def\longra{\longrightarrow}

\font\tenbb=msbm10
\font\sevenbb=msbm7
\font\fivebb=msbm5
\newfam\bbfam
\textfont\bbfam=\tenbb \scriptfont\bbfam=\sevenbb
\scriptscriptfont\bbfam=\fivebb
\def\bb{\fam\bbfam}

\def\Cb{{\bb C}}

\def\Rb{{\bb R}}

\def\Tb{{\bb T}}
\def\Zb{{\bb Z}}

\def\Ac{{\cal A}}

\def\Hc{{\cal H}}

\def\vp{\varphi}

\def\om{\omega}
\def\t{\theta}
\def\d{\delta}

\def\s{\sigma}
\def\G{\Gamma}

\def\g{\gamma}

\def\fl{\forall}

\def\sbs{\subset}

\def\build#1_#2^#3{\mathrel{
\mathop{\kern 0pt#1}\limits_{#2}^{#3}}}

\begin{document}

\title{
\centerline{ NONCOMMUTATIVE MANIFOLDS}
\centerline{\, }
\centerline{THE INSTANTON ALGEBRA}
\centerline{\, }
\centerline{AND ISOSPECTRAL DEFORMATIONS}}

\author{\bf Alain CONNES\refnote{1}, Giovanni LANDI\refnote{2}}

\affiliation{\affnote{1}  Coll\`ege de France, 3, rue Ulm, 75005 PARIS\\
and \\ I.H.E.S., 35, route de Chartres, 91440 BURES-sur-YVETTE \\
\affnote{2} Dipartimento di Scienze Matematiche, Universit\`a di
TRIESTE \\ Via Valerio 12/b, 34127, TRIESTE}

\begin{abstract}
We give new examples of noncommutative manifolds that are less standard
than the
NC-torus or Moyal deformations of $\Rb^n$.
They arise naturally from basic considerations of noncommutative differential
topology and have non-trivial global features.

The new examples include the
instanton algebra  and the NC-4-spheres $S^4_{\theta}$. \\
We construct  the noncommutative algebras $\Ac=C^{\ify} (S^{4}_{\theta})$
of functions on NC-spheres
as solutions to the
vanishing, $ {\rm ch}_j (e) = 0  \, , j < 2 $,
of the Chern character in the cyclic homology of $\Ac$ of an idempotent $e
\in M_4 (\Ac) \, , \, e^2 = e \, , \, e = e^*$.
  We describe the universal noncommutative space obtained from this equation as a
noncommutative Grassmanian as well as the corresponding notion of admissible morphisms.
This space
${\rm Gr}$
contains the suspension of a NC-3-sphere intimately related to
quantum group deformations ${\rm SU}_q (2)$ of ${\rm SU} (2)$ but for
unusual values (complex values of modulus one) of the
  parameter $q$ of $q$-analogues, $q=\exp (2\pi i \t)$.
\noindent

We then construct the noncommutative geometry of $S_{\t}^4$ as
given by a
spectral triple $(\Ac , \Hc , D)$ and check all axioms of noncommutative
manifolds. In a previous paper it was shown that for any Riemannian metric
$g_{\mu \nu}$ on $S^4$ whose volume form $\sqrt{g} \, d^4 x$ is the same as
the one for the round
metric, the corresponding Dirac operator gives a solution to the following
quartic
equation,
\[
\left\langle \left( e - \frac{1}{2} \right) [D,e]^4 \right\rangle = \g_5
\]
where $\langle \ \rangle$ is the projection on the commutant of $4 \ts 4$
matrices.
\\ We shall show how to construct the Dirac operator $D$ on the
noncommutative 4-spheres $S_{\t}^4$ so that the previous equation continues to
hold without any change.
Finally, we show that any compact Riemannian spin manifold whose isometry
group has rank $r \geq 2$ admits
isospectral deformations to noncommutative geometries.
\end{abstract}

\vfill\eject

\section{Introduction}

\noindent It is important to have available examples of noncommutative
manifolds that are less standard than the NC-torus \cite{C3}, \cite{Co-R}
or the
old Moyal deformation of $\Rb^n$ whose algebra is boring.
This is particularly so in view of the upsurge of activity in
the interaction between string theory and  noncommutative geometry started
in \cite{CDS}, \cite{N-S},
\cite{Witten}.
\\ The new examples should arise naturally, have non-trivial global
features (and also
pass the test of noncommutative manifolds as defined in \cite{Co3}). \\
This paper
will provide and analyse very natural such new examples, including the
instanton
algebra and the NC-4-spheres $S^4_{\theta}$, obtained from basic
considerations of
noncommutative differential topology.

\noindent We shall also show quite generally that any compact Riemannian
spin manifold whose
isometry group has rank $r \geq 2$ admits
isospectral deformations to noncommutative geometries.

\noindent  A noncommutative geometry is described by a spectral triple
\begin{equation}
(\Ac , \Hc , D) \label{eq1}
\end{equation}
where $\Ac$ is a noncommutative algebra with involution $ * $, acting
in the
Hilbert
space $\Hc$ while $D$ is a self-adjoint operator with compact resolvent
and such that,
\begin{equation}
[D,a] \ \hbox{is bounded} \ \forall \, a \in \Ac \, . \label{eq2}
\end{equation}
The operator $D$ plays in general the role of the Dirac operator \cite{L-M} in
ordinary Riemannian geometry. It specifies both the metric on the state
space of
$\Ac$ by
\begin{equation}
d (\vp , \psi) = \Sup \, \{ \vert \vp (a) - \psi (a) \vert ; \Vert [D,a] \Vert
\leq 1 \} \label{eq3}
\end{equation}
and the $K$-homology fundamental class (cf. \cite{[Co]}). What
holds things together
in this spectral point of view on NCG is the nontriviality of the pairing
between
the $K$-theory of the algebra $\Ac$ and the $K$-homology class of $D$,
given in
the even case by
\begin{equation}
[e] \in K_0 (\Ac) \ra \hbox{Index} \ D_e^+ \in \Zb \, . \label{eq4}
\end{equation}
Here $[e]$ is the class of an idempotent
\begin{equation}
e \in M_r (\Ac) \, , \ e^2 = e \, , \ e = e^* \label{eq5}
\end{equation}
in the algebra of $r \ts r$ matrices over $\Ac$, and
\begin{equation}
D_e^+ = e \, D^+ e \ , \label{eq6}
\end{equation}
where $ D^+ =
D \, (\frac{1+\g}{2})$
is the restriction of $D$ to the range $\Hc^+$ of $\frac{1+\g}{2}$ and $\g$
is the $\Zb / 2$ grading of $\Hc$ in the even case; thus  $D$ is of the form,
\begin{equation}
D = \left[ \matrix{
0 &D_+^* \cr
D_+ &0 \cr}
\right] \ , \quad \g = \left[ \matrix{ 1 &0 \cr 0&-1 \cr} \right] \, .
\label{eq7}
\end{equation}
The corner stone of the general theory is an operator theoretic index formula
\cite{[Co]}, \cite{[C-M2]}, \cite{Gracia+V} which expresses the above index
pairing (\ref{eq4}) by explicit
{\it local} cyclic cocycles on the algebra $\Ac$. These local formulas become
extremely simple in the special case where only the top component of the Chern
character ${\rm Ch} (e)$ in cyclic homology fails to vanish. (This is easy
to understand in the
analogous simpler case
of ordinary manifolds since the Atiyah-Singer index formula gives the
integral of
the product of the Chern character ${\rm Ch} (E)$, of the bundle $E$ over
the manifold
$M$, by the index class; if the only component of
${\rm Ch} (E)$ is ${\rm ch}_n$, $n = \frac{1}{2} \, \dim M$ only the
$0$-dimensional
component of the index class is involved in the index formula.)

\noindent   Under this assumption the index formula reduces indeed to the
following,
\begin{equation}
\hbox{Index} \, D_e^+ = (-1)^m \int\!\!\!\!\!\!- \ \g \left( e - \frac{1}{2}
\right) [D,e]^{2m} \, D^{-2m} \label{eq8}
\end{equation}
provided the components ${\rm ch}_j (e)$ all vanish for $j < m$. Here $\g$
is the
$\Zb / 2$ grading of $\Hc$ as above, the resolvent of $D$ is of order
$\frac{1}{2m}$ (i.e. its characteristic values $\mu_k$ are
$0(k^{-\frac{1}{2m}})$)
and $\int\!\!\!\!\!-$ is the coefficient of the logarithmic divergency in the
ordinary operator trace \cite{[Dx]} \cite{[Wo]}.

\noindent  We began in \cite{cojmp} to investigate the algebraic relations
implied by the
vanishing,
\begin{equation}
{\rm ch}_j (e) = 0 \qquad j < m \, , \label{eq9bis}
\end{equation}
of the Chern character of $e$ in the cyclic homology of $\Ac$. Note that
this vanishing at the chain level
is a much stronger condition than the vanishing of the usual Chern
differential form.

\noindent  For $m=1$ (and $r = 2$ in (\ref{eq5})) we found commutative
  solutions with $\Ac =
C^{\ify} (S^2)$ as the algebra generated by the matrix components,
\begin{equation}
e_{ij} \, , \ e = [e_{ij}] \in M_2 (\Ac) \, . \label{eq10}
\end{equation}
In fact, for $m=1$ the commutativity is imposed by the relations $e^2 = e$,
$e =
e^*$ and ${\rm ch}_0 (e) = 0$.

\noindent  For $m=2$ (and $r=4$ in (\ref{eq5})) we also found commutative
solutions with $\Ac
= C^{\ify} (S^4)$ where $S^4$ appears as quaternionic projective space but the
computations of \cite{cojmp} used an ``Ansatz'' and did not analyse the
general
solution. In particular this left open the possibility of a noncommutative
solution for $m=2$ (and $r=4$). We shall show in this paper that such
noncommutative
solutions
do exist and provide very natural examples of NC 4-spheres $S_{\t}^4$. We
shall
also describe the noncommutative space associated to (\ref{eq9bis}) for $m=2$
(and $r=4$)
as a
noncommutative Grassmanian as well as the corresponding notion of admissible morphisms.
This space
${\rm Gr}$
contains the suspension of a NC-3-sphere which is intimately related to
quantum group deformations of ${\rm SU} (2)$ but for complex values of
modulus one of the
usual parameter $q$ of $q$-analogues, $q=\exp (2\pi i \t)$.

\noindent  Our next task will be to analyse the metrics (i.e. the operators
$D$) on our
solutions of equation (\ref{eq9bis}).

\noindent  In \cite{cojmp} it was shown that for any Riemannian metric $g_{\mu
\nu}$ on $S^4$ whose
volume form $\sqrt{g} \, d^4 x$ is the same as the one for the round
metric, the
corresponding Dirac operator gives a solution to the following quartic
equation,
\begin{equation}
\left\langle \left( e - \frac{1}{2} \right) [D,e]^4 \right\rangle = \g_5
\label{eq11}
\end{equation}
where $\langle \ \rangle$ is the projection on the commutant of $4 \ts 4$
matrices (recall that $e
\in M_4 (\Ac)$ is a $4 \ts 4$ matrix).

\noindent  We shall show in this paper how to construct the Dirac operator
on the
noncommutative 4-spheres $S_{\t}^4$ so that equation (\ref{eq11}) continues to
hold without any change. Combining this equation (\ref{eq11}) with the index
formula gives a quantization of the volume,
\begin{equation}
\int\!\!\!\!\!\!- \ ds^4 \in \Nb \qquad ds = D^{-1} \label{eq12}
\end{equation}
and fixes (in a given $K$-homology class for the operator $D$) the leading
term of
the spectral
action \cite{[C-C]},
\begin{equation}
\hbox{Trace} \left( f \left( \frac{D}{\Lambda} \right) \right) =
\frac{\Lambda^4}{2} \ \int\!\!\!\!\!\!- \ ds^4 + \cdots \label{eq13}
\end{equation}
Since the next term is the Hilbert-Einstein action in the usual Riemannian
case \cite{[C-C]},
\cite{[Kas]}, \cite{[K-W]},
it is very natural to compare various solutions (commutative or not) of
(\ref{eq11}) using this action.

\vspace{1cm}

\section{Components of the Chern character and the Instanton algebra}

\noindent Let $\Ac$ be an algebra (over $\Cb$) and
\begin{equation}
e \in M_r (\Ac) \, , \ e^2 = e \label{eq2.1}
\end{equation}
be an idempotent.

\noindent The component ${\rm ch}_n (e)$ of the (reduced) Chern character
of $e$ is an
element of
\begin{equation}
\Ac \ot \underbrace{\ov{\Ac} \ot \cdots \ot \ov{\Ac}}_{2n} \label{eq2.2}
\end{equation}
where $\ov{\Ac} = \Ac / \Cb 1$ is the quotient of $\Ac$ by the scalar
multiples of
the unit 1.

\noindent The formula for ${\rm ch}_n (e)$ is (with $\lb_n$ a normalization
constant),
\begin{equation}
{\rm ch}_n (e) = \lb_n \sum \left( e_{i_0 i_1} - \frac{1}{2} \, \d_{i_0 i_1}
\right) \ot e_{i_1 i_2} \ot e_{i_2 i_3} \cdots \ot e_{i_{2n} i_0} \label{eq2.3}
\end{equation}
where $\d_{ij}$ is the usual Kronecker symbol and only the class $\wt{e}_{i_k
i_{k+1}} \in \ov{\Ac}$ is used in the formula. The crucial property of the
components ${\rm ch}_n (e)$ is that they define a {\it cycle} in the $(b,B)$
bicomplex of cyclic homology \cite{Co$_{18}$}, \cite{L},
\begin{equation}
B \, {\rm ch}_n (e) = b \, {\rm ch}_{n+1} (e) \, . \label{eq2.4}
\end{equation}
For any pair of integers $m,r$ we let $\Ac_{m,r}$ be the universal algebra
associated to 
the relations,
\begin{equation}
{\rm ch}_j (e) = 0 \qquad \fl \, j < m \, . \label{eq2.5}
\end{equation}
More precisely we let $\Ac_{m,r}$ be generated by the $r^2$ elements $e_{ij}$;
$i,j \in \{ 1, \ldots , r \}$ and we first impose the relations
\begin{equation}
e^2 = e \qquad e = [e_{ij}] . \label{eq2.6}
\end{equation}
An admissible homomorphism,
\begin{equation}
\rho : \Ac_{m,r} \ra \Bc, \label{eq2.7}
\end{equation}
to an arbitrary algebra $\Bc$, is
given
by the $\rho (e_{ij}) \in \Bc$ which fulfill
\begin{equation}
\rho (e)^2 = \rho (e) \, , \label{eq2.8}
\end{equation}
and ${\rm ch}_j (\rho(e)) = 0 $ for $ \, j < m$, thus 
\begin{equation}
\sum \left( \rho (e_{i_0 i_1}) - \frac{1}{2} \ \d_{i_0 i_1} \right) \ot
\wt{\rho
(e_{i_1 i_2})} \ot \cdots \ot \wt{\rho (e_{i_{2j} i_0})} = 0 \label{eq2.9}
\end{equation}
where the symbol $\, ^\sim \, $ means that only the class in $\ov{\Bc}$
matters.
We define
$\Ac_{m,r}$
as the quotient of the algebra defined by (\ref{eq2.6}) by the intersection of
kernels of all admissible morphisms $\rho$.

\noindent Elements of the algebra $\Ac_{m,r}$ can be represented as
polynomials in the
generators
$e_{ij}$ and to prove that such a polynomial $P (e_{ij})$ is non zero in
$\Ac_{m,r}$ one must
construct a solution to the
above equations for which $P (e_{ij}) \ne 0$.

\noindent To get a $C^*$ algebra we endow $\Ac_{m,r}$ with the involution
given by,
\begin{equation}
(e_{ij})^* = e_{ji} \label{eq2.10}
\end{equation}
which means that $e = e^*$ in $M_r (\Ac)$. We define a norm by,
\begin{equation}
\Vert P \Vert = \Sup \, \Vert (\pi (P)) \Vert \label{eq2.11}
\end{equation}
where $\pi$ ranges through all representations of the above equations in
Hilbert
space. Such a $\pi$ is given by a Hilbert space $\Hc$ and a self-adjoint
idempotent,
\begin{equation}
E \in M_r (\Lc (\Hc)) \, , \ E^2 = E \, , \ E = E^* \label{eq2.12}
\end{equation}
such that (\ref{eq2.9}) holds for $\Bc = \Lc (\Hc)$.

\noindent One checks that for any polynomial $P (e_{ij})$ the quantity
(\ref{eq2.11}), i.e.
the supremum of the norms,
\begin{equation}
\Vert P (E_{ij}) \Vert \label{eq2.13}
\end{equation}
is finite.

\noindent We let $A_{m,r}$ be the $C^*$ algebra obtained as the completion
of $\Ac_{m,r}$
for the above norm.

\noindent To get familiar with the (a priori noncommutative) spaces ${\rm
Gr}_{m,r}$ such
that,
\begin{equation}
A_{m,r} = C ({\rm Gr}_{m,r}) \label{eq2.14}
\end{equation}
we shall first recall from \cite{cojmp} what happens in the simplest case
$m=1$, $r=2$.

\noindent One has $e = \left[ \matrix{ e_{11} &e_{12} \cr e_{21} &e_{22}
\cr} \right]$ and
the condition (\ref{eq2.7}) just means that
\begin{equation}
e_{11} + e_{22} = 1 \label{eq2.15}
\end{equation}
while (\ref{eq2.6}) means that
\begin{eqnarray}
&e_{11}^2 + e_{12} \, e_{21} = e_{11} \, , \ e_{11} \, e_{12} + e_{12} \,
e_{22} =
e_{12} \, , \label{eq2.16} \\
&e_{21} \, e_{11} + e_{22} \, e_{21} = e_{21} \, , \ e_{21} \, e_{12} +
e_{22}^2 =
e_{22} \, . \nonumber
\end{eqnarray}

\noindent By (\ref{eq2.15}) we get $e_{11} - e_{11}^2 = e_{22} - e_{22}^2$,
so that
(\ref{eq2.16}) shows that $e_{12} \, e_{21} = e_{21} \, e_{12}$. We also
see that
$e_{12}$ and $e_{21}$ both commute with $e_{11}$. This shows that
$\Ac_{1,2}$ is
commutative and allows to check that ${\rm Gr}_{1,2} = S^2$ is the
2-sphere. Thus
${\rm Gr}_{1,2}$ is an ordinary commutative space.

\noindent
Next, we move on to the case $m=2$, $r=4$. Our main task now will be to
show that ${\rm Gr}_{2,4}$ is a very interesting noncommutative space.
Note that the notion of admissible 
morphism is a non trivial piece of structure on ${\rm Gr}_{2,4}$
since the identity map is not admissible.

\noindent We can first reformulate the construction of \cite{cojmp} section
XI and get an admissible 
surjection,
\begin{equation}
A_{2,4} \build \longra_{}^{\s} C (S^4) \label{eq2.17}
\end{equation}
where $S^4$ appears naturally as quaternionic projective space, $S^4 = P_1
(\Hb)$.

\noindent Let us recall from \cite{cojmp} that the equality,
\begin{equation}
E(x) = \left[ \matrix{
t &q \cr
\ov q &1-t \cr
} \right] \in M_4 (\Cb) \label{eq2.18}
\end{equation}
for $x = (q,t)$ given by a pair of a quaternion
$q = {\scriptstyle        
  \addtolength{\arraycolsep}{-.5\arraycolsep}
  \renewcommand{\arraystretch}{0.5}
  \left[ \begin{array}{cc}
  \scriptstyle \a  & \scriptstyle \b \\
  \scriptstyle -\b^*  & \scriptstyle \a^*  \end{array} \scriptstyle\right] }$
and a real number $t$ such that
\begin{equation}
q \, \ov q = t - t^2 \label{eq2.19}
\end{equation}
defines a map from the 4-sphere $S^4$ (the double of the 4-disk $\vert q \vert
\leq 1$) to the Grassmanian of 2-dimensional projections $E = E^2 = E^*$ in
$M_4 (\Cb)$ such that,
\begin{equation}
\hbox{Trace} \, (F(x) \, F(y) \, F(z)) = 0 \qquad \fl \, x,y,z \in S^4
\label{eq2.20}
\end{equation}
where $F(x) = 2 \, E(x) - 1$ is the corresponding self-adjoint isometry.

\noindent The equality XI.54 of \cite{cojmp} is weaker than this statement
but examining the
proof one gets (\ref{eq2.20}). To formulate the result for arbitrary even
spheres
$S^{2m}$ we note first that using (\ref{eq2.4}) the equality
\begin{equation}
\om = {\rm ch}_m (e) \label{eq2.21}
\end{equation}
defines a {\it Hochschild cycle} $ \rho(\om) \in Z_{2m} (\Bc)$ for any 
admissible morphism $\rho : \Ac_{m,r} \ra \Bc$. We let $r =
2^m$ and
construct an admissible surjection,
\begin{equation}
A_{m,2^m} \build \longrightarrow_{}^{\s} C (S^{2m}) \label{eq2.22}
\end{equation}
which is non trivial inasmuch as
\begin{equation}
\s (\om) = v \label{eq2.23}
\end{equation}
is the volume form of the round oriented sphere.

\noindent To construct $\s$ we let $C = {\rm Cliff} \, (\Rb^{2m})$ be the
Clifford algebra
of the (oriented) Euclidean space $\Rb^{2m}$. We identify $S^{2m}$ with the
space
of pairs $(\xi , t)$, $\xi \in \Rb^{2m}$ and $t \in [-1,1]$ such that
$\Vert \xi
\Vert^2 + t^2 = 1$. We then define a map from $S^{2m}$ to the Grassmanian of
self-adjoint idempotents in $C$ by
\begin{equation}
E (\xi , t) = \frac{1}{2} + \frac{1}{2} \, (\g (\xi) + t \, \g) \label{eq2.24}
\end{equation}
where $\g (\xi)$ is the usual inclusion of $\Rb^{2m}$ in $C$ such that
\begin{equation}
\g (\xi)^2 = \Vert \xi \Vert^2 \, , \ \g (\xi) = \g (\xi)^* \label{eq2.25}
\end{equation}
and $\g \in C$, $\g^* = \g$, $\g^2 = 1$ is the $\Zb / 2$ grading
associated
with the chosen orientation of $\Rb^{2m}$. One has $\g \, \g (\xi) = -\g
(\xi) \g$
for 	any $\xi$ which allows to check that $\g (\xi) + t \, \g$ is an
involution
and $E$ a self-adjoint idempotent. Next, for $\ell < 2m$, $\ell$ {\it odd},
\begin{equation}
\hbox{Trace} \, ((\g (\xi_1) + t_1 \, \g) \ldots (\g (\xi_{\ell}) +
t_{\ell} \,
\g)) = 0 \qquad \fl \, \xi_j , t_j \, . \label{eq2.26}
\end{equation}
Indeed the coefficient of monomials in $t$ of even degree is of the form
$\hbox{Trace} \, (\g (\xi_1) \ldots$ $\g (\xi_{2k+1}))$ which vanishes by
anticommutation with $\g$. The coefficient of monomials in $t$ of odd
degree is of
the form $\hbox{Trace} \, (\g (\xi_1) \ldots \g (\xi_{2k}) \, \g)$ where $k
< m$.
It vanishes because $\g$ is orthogonal to all the lower filtration of $C$.
We thus
get,
\begin{equation}
\hbox{Trace} \, \left(\left(E (x_1) - \frac{1}{2} \right) \ldots \left( E
(x_{\ell}) - \frac{1}{2} \right) \right) = 0 \qquad \fl \, x_1 , \ldots ,
x_{\ell}
\in S^{2m} \label{eq2.27}
\end{equation}
provided $\ell$ is odd, $\ell < 2m$.

\noindent Hence $E$ defines an admissible homomorphism $\s : A_{m,2^m} \ra C(S^{2m})$
and one has, as in
\cite{cojmp}, the following result,
\bigskip

\begin{theorem}\label{th1}
\begin{itemize}
\item[{\rm a)}] $E \in C^{\ify} (S^{2m} , M_r (\Cb))$ satisfies $E = E^2 =
E^*$
and ${\rm ch}_j (E) = 0$ $\fl \, j < m$.
\item[{\rm b)}] The Hochschild cycle $\om = {\rm ch}_m (E)$ is the volume
form of
the round sphere $S^{2m}$.
\item[{\rm c)}] Let $g$ be a Riemannian metric on $S^{2m}$ with volume form
$\sqrt{g} \, d^{2m} x = \om$, then the corresponding Dirac operator $D$
fulfills
$$
\left\langle \left( e - \frac{1}{2} \right) [D,e]^{2m} \right\rangle = \g
$$
where $e=E$ as above and $\langle \ \rangle$ is the projection on the
commutant of
$M_r (\Cb)$.
\end{itemize}
\end{theorem}

\noindent We have identified $M_r (\Cb)$ with the Clifford algebra $C$, $r
= 2^m$. This
result shows in particular that ${\rm Gr}_{m,r}$, $r=2^m$, contains $S^{2m}$ in
such a way that $\om | S^{2m}$ is the volume form for the round metric. The
proof
is the same as in \cite{cojmp}.

\vspace{1cm}

\section{The noncommutative 4-sphere}

\noindent Let us now move on to the inclusion $S_{\t}^4 \sbs {\rm
Gr}_{2,4}$ where
$S_{\t}^4$ is the noncommutative 4-sphere we are about to describe.

\noindent One should observe from the outset that the compact Lie group
$SU(4)$ acts by
automorphisms,
\begin{equation}
PSU(4) \subset {\rm Aut} \, (C^{\infty} {\rm Gr}_{2,4}) \label{eq9}
\end{equation}
by the following operation,
\begin{equation}
e \ra U \, e \, U^*
\end{equation}
where $U \in SU(4)$ is viewed as a $4 \times 4$ matrix and $e = [e_{ij}]$ as
above.

\smallskip

\noindent  We shall now show that the algebra $C({\rm Gr_{2,4}})$ is
noncommutative by constructing explicit admissible surjections,
\begin{equation}
C ({\rm Gr_{2,4}}) \ra C(S_{\theta}^4)
\end{equation}
whose form is dictated by natural deformations of the 4-sphere similar in
spirit to the standard deformation of $\Tb^2$ to $\Tb_{\theta}^2$.

\noindent We first determine the algebra generated by $M_4 (\Cb)$ and a
projection $e = e^*$ such that $\left\lgl e - \frac{1}{2} \right\rgl = 0$ as
above and whose two by two matrix expression is of the form,
\begin{equation}
[ e^{ij} ] = \left[ \matrix{q_{11} &q_{12} \cr q_{21} &q_{22} \cr} \right]
\label{eq58}
\end{equation}
where each $q_{ij}$ is a $2 \ts 2$ matrix of the form,
\begin{equation}
q = \left[ \matrix{\a &\b \cr -\lambda \b^* &\a^* \cr} \right] \, ,
\label{eq59}
\end{equation}
and $\lambda=\exp (2\pi i \t)$ is a complex number of modulus one,
different from -1 for convenience.
Since $e = e^*$, both $q_{11}$ and $q_{22}$ are selfadjoint, moreover since
$\left\lgl e - \frac{1}{2} \right\rgl = 0$, we can
find $t = t^*$ such that,
\begin{equation}
q_{11} = \left[ \matrix{t &0 \cr 0 &t \cr} \right] \, , \ q_{22} = \left[
\matrix{(1-t) &0 \cr 0 &(1-t) \cr} \right] \, . \label{eq60}
\end{equation}
We let
$q_{12} = {\scriptstyle        
  \addtolength{\arraycolsep}{-.5\arraycolsep}
  \renewcommand{\arraystretch}{0.5}
  \left[ \begin{array}{cc}
  \scriptstyle \a  & \scriptstyle \b \\
  \scriptstyle -\lambda \b^*  & \scriptstyle \a^*  \end{array}
\scriptstyle\right] }$ ,
we then get from $e = e^*$,
\begin{equation}
q_{21} = \left[ \matrix{\a^* &- \bar{\lambda} \b \cr \b^* &\a \cr} \right]
\, . \label{eq61}
\end{equation}
We thus see that the commutant $\Bc_{\theta}$ of $M_4 (\Cb)$ is generated
by $t,\a,\b$
and we first need to find the relations imposed by the equality $e^2 = e$.

\noindent In terms of
\begin{equation}
e = \left[ \matrix{ t &q \cr q^* &1-t \cr} \right] \, , \label{idempot}
\end{equation}
the equation $e^2 = e$ means that $t^2 - t + q q^* = 0$, $t^2 - t + q^* q =
0$ and $[t,q] = 0$. This shows that $t$ commutes with $\a$, $\b$, $\a^*$ and
$\b^*$ and since $qq^* = q^* q$ is a diagonal matrix
\begin{equation}
\a \a^* = \a^* \a \, , \ \a \b =\lambda  \b \a \, , \ \a^* \b
=\bar{\lambda} \b \a^* \, , \ \b
\b^* = \b^* \b \label{eq62}
\end{equation}
so that the $C^*$ algebra $\Bc_{\theta}$ is not commutative for $ \lambda$
different from
1. The only further relation is,
(besides $t = t^*$),
\begin{equation}
\a \a^* + \b \b^* + t^2 - t = 0 \, . \label{eq63}
\end{equation}
We denote by $S^4_{\theta}$ the corresponding noncommutative space,
so that $C(S^4_{\theta})=\Bc_{\theta}$.
It is by construction the suspension of the noncommutative 3-sphere
$S^3_{\theta}$
whose coordinate algebra is
generated by $\a$ and $\b$ as above and say the special value $t=1/2$.
This noncommutative 3-sphere is related by analytic continuation of the
parameter
to the quantum
group $SU(2)_q$ but the usual theory requires $q$ to be real whereas we
need a complex
number of
modulus one which spoils the unitarity of the coproduct \cite{wor}.
\\
Had we taken the deformation parameter to be real, $\lambda=q\in
\Rb$, like in
\cite{dlm} we would have obtained a different deformation $S^4_{q}$ of the
commutative
sphere $S^4$, whose algebra is different from the above one. More
important, the two
dimensional component $ ch_1(e)$ of the Chern character
would not
vanish.

\bigskip
\noindent We shall now check that for the sphere $S^4_{\theta}$ the two
dimensional
component $ ch_1(e)
$ automatically vanishes as an element of the (normalized)
(b,B)-bicomplex so that,
\begin{equation}
  ch_{n}(e)  = 0 \, , \
n = 0,1  . \label{eq67}
\end{equation}
With
$q = {\scriptstyle        
  \addtolength{\arraycolsep}{-.5\arraycolsep}
  \renewcommand{\arraystretch}{0.5}
  \left[ \begin{array}{cc}
  \scriptstyle \a  & \scriptstyle \b \\
  \scriptstyle -\lambda \b^*  & \scriptstyle \a^*  \end{array}
\scriptstyle\right] }$ ,
we get,
\begin{eqnarray}
& ch_1(e)  &= \Biggl\lgl \left( t - \frac{1}{2} \right)
\, (dq \, dq^* - dq^* \, dq) \\
&&+ \, q \, (dq^* \, dt - dt \, dq^*) + q^* \, (dt \, dq - dq \, dt )
\Biggl\rgl \nonumber
\label{eq69}
\end{eqnarray}
where the expectation in the right hand side is relative to $M_2 (\Cb)$ and
we use
the notation $\, d \, $ instead of the tensor notation.

\noindent The diagonal elements of $\om = dq \, dq^*$ are
$$
\om_{11} = d \a \, d\a^* + d\b \, d\b^* \, , \ \om_{22} = d\b^* \, d\b +
d\a^* \, d\a
$$
while for $\om' = dq^* \, dq$ we get,
$$
\om'_{11} = d\a^* \, d\a + d\b \, d\b^* \, , \ \om'_{22} = d\b^* \, d\b + d\a
\, d\a^* \, .
$$
It follows that, since $t$ is diagonal,
\begin{equation}
\left\lgl \left( t - \frac{1}{2} \right) \, (dq \, dq^* - dq^* \, dq)
\right\rgl = 0 \, . \label{eq70}
\end{equation}
The diagonal elements of $q \, dq^* \, dt = \rho$ are
$$
\rho_{11} = \a \, d\a^* \, dt + \b \, d\b^* \, dt \, , \ \rho_{22} = \b^* \,
d\b \, dt + \a^* \, d\a \, dt
$$
while for $\rho' = q^* \, dq \, dt$ they are
$$
\rho'_{11} = \a^* \, d\a \, dt + \b \, d\b^* \, dt \, , \ \rho'_{22} = \b^* \,
d\b \, dt + \a \, d\a^* \, dt \, .
$$
Similarly for $\s = q \, dt \, dq^*$ and $\s' = q^* \, dt \, dq$ one gets the
required cancellations so that,
\begin{equation}
  ch_1(e)  = 0 \, ,
\label{eq71}
\end{equation}
We thus get,
\bigskip

\begin{theorem}\label{th2}
\begin{itemize}
\item[{\rm a)}] $e \in C^{\ify} (S^4_{\theta}, M_4 (\Cb))$ satisfies $e =
e^2 = e^*$
and ${\rm ch}_j (e) = 0$ $\fl \, j < 2$.
\item[{\rm b)}] ${\rm Gr}_{2,4}$ is a noncommutative space and
$S^4_{\theta} \sbs {\rm Gr}_{2,4} $.

\end{itemize}
\end{theorem}

\vspace{1cm}

\noindent Since $ ch_1(e)=0$, it follows that $ ch_2(e)$ is a Hochschild
cycle which will play the role of the
round volume form on $S^4_{\theta}$ and that we shall now compute.
With the above notations one has,
\begin{equation}
  ch_2(e)
  = \left\lgl  \left[ \matrix{ t - \frac{1}{2} &q \cr q^* &\frac{1}{2}-t \cr}
\right]
(\left[ \matrix{ dt  &dq \cr dq^* & -dt \cr} \right])^4\right\rgl
\label{ch4}
\end{equation}
and the sum of the diagonal elements is
\begin{eqnarray}
&& \left( t - \frac{1}{2} \right) \, \Big( (dt^2 + dq \, dq^*)^2 + (dt \,
dq - dq
\, dt) (dq^* \, dt - dt \, dq^*) \Big) \label{ch4.1}
\\
- && \left( t - \frac{1}{2} \right) \, \Big( (dq^* \, dt - dt \, dq^*) (dt
\, dq -
dq \, dt) + (dq^* \, dq + dt^2)^2 \Big) \nonumber \\
+ && q \Big( (dq^* \, dt - dt \, dq^*) (dt^2 + dq \, dq^*) + (dq^* \, dq +
dt^2) (dq^* \, dt - dt \, dq^*) \Big) \nonumber \\
+ && q^* \Big((dt^2 + dq \, dq^*) (dt \, dq - dq \, dt) + (dt \, dq - dq \, dt)
(dq^* \, dq + dt^2) \Big) \, . \nonumber
\end{eqnarray}
Since $t$ and $dt$ are diagonal $2\times 2$ matrices of operators and the same
diagonal terms appear in $dq \, dq^*$ and $dq^* \, dq$, by the same argument by
which we got the vanishing (\ref{eq70}), the first two lines only
contribute by,
\begin{equation}
\left\lgl \left( t - \frac{1}{2} \right) \, (dq \, dq^* \, dq \, dq^* -
dq^* \, dq
\, dq^* \, dq)
\right\rgl \, . \label{ch4.1first}
\end{equation}
Similarly, the last two lines only contribute by
\begin{eqnarray}
&& \Big\lgl q^* \, (dt \, dq \, dq^* \, dq  - dq \, dt \, dq^* \, dq +
dq \, dq^* \, dt \, dq - dq \, dq^* \, dq \, dt) \label{ch4.1last}
\\ && - q \, (dt \, dq^* \, dq \, dq^* - dq^* \, dt \, dq \, dq^* + dq^* \,
dq \,
dt \, dq^* - dq^* \, dq \, dq^* \, dt)
\Big\rgl \, . \nonumber
\end{eqnarray}
The direct computation gives $ ch_2(e) $ as a sum of five components
\begin{equation}
  ch_2(e)  = (t - \frac{1}{2}) \, \G_t + \a \, \G_{\a} + \a^* \, \G_{\a^*} +
\b \,
\G_{\b} +
\b^* \, \G_{\b^*} \, \label{gamma5},
\end{equation}
with the operators $\G_t \, , \G_{\a} \, , \G_{\a^*} \, , \G_{\b} \, ,
\G_{\b^*}$
explicitly given by

\begin{eqnarray}
\G_t &=&
(d\a \, d\a^* - d\a^* \, d\a)(d\b \, d\b^* - d\b^* \, d\b) \label{gammat}  \\
&\, & \, \, \, \, \, \, \, \, \, \,
+ \, (d\b \, d\b^* - d\b^* \, d\b) (d\a \, d\a^* - d\a^* \, d\a) \nonumber \\
&+& \, (d\a \, d\b - \lambda \, d\b \, d\a)
(d\b^* \, d\a^* - \bar{\lambda} \, d\a^* \, d\b^*) \nonumber \\
&\, & \, \, \, \, \, \, \, \, \, \,
+ \, (d\b^* \, d\a^* - \bar{\lambda} \, d\a^* \, d\b^*)
(d\a \, d\b - \lambda \, d\b \, d\a) \nonumber \\
&+& \, (d\a^* \, d\b - \bar{\lambda} \, d\b \, d\a^*)
(\lambda \, d\a \, d\b^* \, - d\b^* \, d\a) \nonumber \\
&\, & \, \, \, \, \, \, \, \, \, \,
+ \, (\lambda \, d\a \, d\b^* \, -  d\b^* \, d\a)
(d\a^* \, d\b - \bar{\lambda} \, d\b \, d\a^*) \, ; \nonumber
\end{eqnarray}

\begin{eqnarray}
\G_\a &=&
(dt \, d\a^* - d\a^* \, dt) (d\b^* \, d\b - d\b \, d\b^*) \label{gammaa}  \\
&\, & \, \, \, \, \, \, \, \, \, \,
+ \, (d\b^* \, d\b - d\b \, d\b^*) (dt \, d\a^* - d\a^* \, dt) \nonumber \\
&+& \, (d\b \, dt - dt \, d\b)
(d\b^* \, d\a^* - \bar{\lambda} \, d\a^* \, d\b^*) \nonumber \\
&\, & \, \, \, \, \, \, \, \, \, \,
+ \, \lambda \, (d\b^* \, d\a^* - \bar{\lambda} \, d\a^* \, d\b^*)
(d\b \, dt - dt \, d\b) \nonumber \\
&+& \, (d\a^* \, d\b - \bar{\lambda} \, d\b \, d\a^*)
(d\b^* \, dt - dt \, d\b^*) \nonumber \\
&\, & \, \, \, \, \, \, \, \, \, \,
+ \, \lambda \, (d\b^* \, dt - dt \, d\b^*)
(d\a^* \, d\b - \bar{\lambda} \, d\b \, d\a^*) \, ; \nonumber
\end{eqnarray}

\begin{eqnarray}
\G_{\a^*} &=&
(dt \, d\a - d\a \, dt) (d\b \, d\b^* - d\b^* \, d\b) \label{gammaabar} \\
&\, & \, \, \, \, \, \, \, \, \, \,
+ \, (d\b \, d\b^* - d\b^* \, d\b) (dt \, d\a - d\a \, dt) \nonumber \\
&+& \, (d\a \, d\b - \lambda \, d\b \, d\a)
(dt \, d\b^* - d\b^* \, dt) \nonumber \\
&\, & \, \, \, \, \, \, \, \, \, \,
+ \, \bar{\lambda} \, (dt \, d\b^* - d\b^* \, dt)
(d\a \, d\b - \lambda \, d\b \, d\a) \nonumber \\
&+& \, (dt \, d\b - d\b \, dt)
(d\b^* \, d\a - \lambda \, d\a \, d\b^*) \,
\nonumber \\
&\, & \, \, \, \, \, \, \, \, \, \,
+ \, \bar{\lambda} \, (d\b^* \, d\a - \lambda \, d\a \, d\b^*)
(dt \, d\b - d\b \, dt) ; \nonumber
\end{eqnarray}

\begin{eqnarray}
\G_\b &=&
(dt \, d\b^* - d\b^* \, dt) (d\a^* \, d\a - d\a \, d\a^*) \label{gammab}  \\
&\, & \, \, \, \, \, \, \, \, \, \,
+ \, (d\a^* \, d\a - d\a \, d\a^*) (dt \, d\b^* - d\b^* \, dt) \nonumber \\
&+& \, \lambda \, (dt \, d\a - d\a \, dt)
(d\b^* \, d\a^* - \bar{\lambda} \, d\a^* \, d\b^*)
\nonumber \\
&\, & \, \, \, \, \, \, \, \, \, \,
+ \, (d\b^* \, d\a^* - \bar{\lambda} \, d\a^* \, d\b^*)
(dt \, d\a - d\a \, dt) \nonumber \\
&+& \, \bar{\lambda} \, (d\a^* \, dt - dt \, d\a^*)
(d\b^* \, d\a - \lambda  \, d\a \, d\b^*) \nonumber \\
&\, & \, \, \, \, \, \, \, \, \, \,
+ \, (d\b^* \, d\a - \lambda  \, d\a \, d\b^*)
(d\a^* \, dt - dt \, d\a^*) \, ; \nonumber
\end{eqnarray}

\begin{eqnarray}
\G_{\b^*} &=&
(dt \, d\b - d\b \, dt) (d\a \, d\a^* - d\a^* \, d\a) \label{gammabbar} \\
&\, & \, \, \, \, \, \, \, \, \, \,
+ \, (d\a \, d\a^* - d\a^* \, d\a) (dt \, d\b - d\b \, dt) \nonumber \\
&+& \, (d\a^* \, dt - dt \, d\a^*)
(d\a \, d\b - \lambda d\b \, d\a)
\nonumber \\
&\, & \, \, \, \, \, \, \, \, \, \,
+ \, \bar{\lambda} \, (d\a \, d\b - \lambda d\b \, d\a)
(d\a^* \, dt - dt \, d\a^*) \nonumber \\
&+& \, (dt \, d\a - d\a \, dt)
(d\a^* \, d\b - \bar{\lambda} \, d\b \, d\a^*) \nonumber \\
&\, & \, \, \, \, \, \, \, \, \, \,
+ \, \lambda \, (d\a^* \, d\b - \bar{\lambda} \, d\b \, d\a^*)
(dt \, d\a - d\a \, dt) \, .  \nonumber
\end{eqnarray}

\noindent One can equivalently (in order to avoid any confusion with
ordinary differentials)
write the Hochschild cycle $c= ch_2(e)$
as
\begin{equation}
c = (t - \frac{1}{2}) \, c_t + \a \, c_{\a} + \a^* \, c_{\a^*} + \b \, c_{\b} +
\b^* \, c_{\b^*} \, ; \label{hocy}
\end{equation}
where the components $c_t \, , c_{\a} \, , c_{\a^*} \, , c_{\b} \, ,
c_{\b^*}$, which are
elements in $\Bc_{\theta}\ot\Bc_{\theta}\ot\Bc_{\theta}\ot\Bc_{\theta}$,
have an
expression of the same form as the corresponding operators in
(\ref{gammat}-\ref{gammabbar}) with the symbol $\, d \, $ substituted by
the tensor
product symbol $\, \ot \, $. The vanishing of $b c$, which has six hundred
terms,
can be checked directly from the commutation relations (\ref{eq62}). The
cycle $c$ is totally
`$\lambda$-antisymmetric'.

\vspace{1cm}

\section{The noncommutative geometry of $S^4_{\theta}$}

\noindent The next step consists in finding the Dirac operator
which gives a solution to the basic quartic equation (I.11).
Let $\Ac = C^{\ify} (S_{\t}^4)$ be the algebra of smooth functions on the
noncommutative sphere $S_{\t}^4$. We shall now describe a spectral triple
\begin{equation}
(\Ac , \Hc , D) \label{eq1bis}
\end{equation}
which describes the geometry on $S_{\t}^4$ corresponding to the round metric.

\noindent In order to do that we first need to find good coordinates on
$S_{\t}^4$ in terms of
which the operator $D$ will be easily expressed. We choose to parametrize
$\a , \b$ and
$t$ as follows,
\begin{equation}
\a = \,  \frac{u}{2} \, \cos \vp \cos \psi  \ , \ \b = \,  \frac{v}{2} \,
\sin \vp \cos
\psi  \, , \ t = \frac{1}{2} + \frac{1}{2} \, \sin \psi \, . \label{eq2bis}
\end{equation}
Here $\vp$ and $\psi$ are ordinary angles with domain
\begin{equation}
0 \leq \vp \leq \frac{\pi}{2} \, , \ - \frac{\pi}{2} \leq \psi \leq
\frac{\pi}{2} \, ,
\label{eq3bis}
\end{equation}
while $u$ and $v$ are the usual unitary generators of the algebra $C^{\ify}
(\Tb_{\theta}^2)$ of smooth functions on the noncommutative 2-torus. Thus
the presentation
of their relations is
\begin{equation}
uv = \lb v u \, , \ uu^* = u^* u = 1 \, , \ vv^* = v^* v = 1 \, .
\label{eq4bis}
\end{equation}
One checks that $\a , \b , t$ given by (2) satisfy the basic presentation
of the
generators of $C^{\ify} (S_{\t}^4)$ which thus appears as a {\it
subalgebra} of the
algebra generated (and then closed under smooth calculus) by $e^{i\vp}$,
$e^{i\psi}$,
$u$ and $v$.

\noindent For $\t = 0$ one readily computes the round metric,
\begin{equation}
G = 4 \, (d\a \, d \ov \a + d \b \, d \ov \b + dt^2) \label{eq5bis}
\end{equation}
and in terms of the coordinates, $\vp , \psi , u , v$ one gets,
\begin{equation}
G = \cos^2 \vp \cos^2 \psi \, du \, d \ov u + \sin^2 \vp \cos^2 \psi \, dv
\, d \ov v +
\cos^2 \psi \, d \vp^2 + d \psi^2 \, . \label{eq6bis}
\end{equation}
Up to normalization, its volume form is given by
\begin{equation}
\sin \vp \cos \vp \, (\cos \psi)^3 \, \ov u \, du \wedge \ov v \, dv \wedge
d\psi \wedge d\vp \, . \label{eq7bis}
\end{equation}
In terms of these rectangular coordinates we get the following simple
expression for
the Dirac operator,
\begin{eqnarray}
D &= &(\cos \vp \cos \psi)^{-1} \, u \, \frac{\partial}{\partial u} \, \g_1
+ (\sin \vp
\cos \psi)^{-1} \, v \, \frac{\partial}{\partial v} \, \g_2 +  \\
&+ &\frac{1}{\cos \psi} \, \sqrt{-1} \, (\frac{\partial}{\partial \vp} \, +
\frac{1}{2} \, {\rm
cotg} \, \vp - \frac{1}{2} \,
{\rm tg} \, \vp) \,
\g_3
+
\sqrt{-1} \, (\frac{\partial}{\partial \psi} \, - \, \frac{3}{2} \, {\rm tg} \,
\psi) \, \g_4 \, . \nonumber
\end{eqnarray}
Here $\g_{\mu}$ are the usual Dirac $4 \ts 4$ matrices with
\begin{equation}
\{ \g_{\mu} , \g_{\nu} \} = 2 \, \d_{\mu \nu} \, , \ \g_{\mu}^* = \g_{\mu}
\, .
\label{eq9ter}
\end{equation}

\noindent It is now easy to move on to the noncommutative case, the only
tricky point
is that there are nontrivial boundary conditions for the operator $D$, but
we shall just
leave them unchanged in the NC case, the only thing which changes is the
algebra and the way it acts in
the Hilbert space.
The formula for the operator $D$ is now,
\begin{eqnarray}
D &= &(\cos \vp \cos \psi)^{-1} \,  \delta_1 \, \g_1 + (\sin \vp
\cos \psi)^{-1} \, \delta_2  \, \g_2 + \label{eq8bis} \\
&+ &\frac{1}{\cos \psi} \, \sqrt{-1} \, (\frac{\partial}{\partial \vp} \, +
\frac{1}{2} \, {\rm
cotg} \, \vp - \frac{1}{2} \,
{\rm tg} \, \vp) \,
\g_3  + \sqrt{-1} \, (\frac{\partial}{\partial \psi} \, - \, \frac{3}{2} \,
{\rm tg} \,
\psi) \, \g_4 \, . \nonumber
\end{eqnarray}
where the $\gamma_\mu$ are the usual Dirac matrices and where $\delta_1$ and
$\delta_2$ are the standard derivations of the NC torus so that
\begin{eqnarray}
&& \delta_1(u) = u \, , \, \, \delta_1(v) = 0 \, ,  \\
&& \delta_2(u) = 0 \, , \, \, \delta_2(v) = v \, ; \nonumber
\end{eqnarray}
One can then check that the
corresponding metric is the right one (the round
metric).

\noindent In order to compute the operator $\left\langle \left( e -
\frac{1}{2} \right) [D,e]^{4} \right\rangle$
(in the tensor product by $M_4( \Cb)$) we need
the commutators of $D$ with the generators of $C^{\ify} (S^{4}_{\theta})$.
They are
given by the
following simple expressions,
\begin{eqnarray}
&&[ D, \a ] = \frac{u}{2} \, \{\gamma_1 -
\, \sqrt{-1} \, \sin(\phi) \, \gamma_3 - \, \sqrt{-1} \, \cos\phi)
\sin(\psi) \, \gamma_4
\, \} \, , \\
&& [ D, \a^* ] = - \frac{u^*}{2} \, \{\gamma_1 +
\, \sqrt{-1} \, \sin(\phi) \, \gamma_3 + \, \sqrt{-1} \, \cos(\phi)
\sin(\psi) \, \gamma_4
\, \} \, ,  \nonumber \\
&& [ D, \b ] = \frac{v}{2} \,  \{ \, \gamma_2
+ \, \sqrt{-1} \, \cos(\phi)  \, \gamma_3 - \, \sqrt{-1} \, \sin(\phi)
\sin(\psi) \, \gamma_4
\, \} \, ,
\nonumber \\
&&  [ D, \b^* ] = - \frac{v^* }{2} \, \{ \,
\gamma_2 - \, \sqrt{-1} \, \cos(\phi)  \, \gamma_3 + \, \sqrt{-1} \,
\sin(\phi) \sin(\psi)  \,
\gamma_4 \,
\} \, , \nonumber \\
&& [ D, t ] = \frac{\, \sqrt{-1} }{2} \cos(\psi) \, \gamma_4 \, . \nonumber
\end{eqnarray}

\noindent We check in particular that they are all bounded operators and
hence that
for any $f \in C^{\ify} (S^{4}_{\theta})$ the commutator $[D, \, f]$ is
bounded.
Then, a long but straightforward calculation shows that the operator
$\left\langle \left( e - \frac{1}{2} \right) [D,e]^{4} \right\rangle$ is
a multiple of $\gamma_5 = \gamma_1 \gamma_2 \gamma_3 \gamma_4$.
One first checks that it is equal to $\pi (c)$ where $c$ is the Hochschild
cycle in (III.\ref{hocy})
  and $\pi$ is the canonical map from the
Hochschild chains to operators given by
\begin{equation}
\pi( a_0 \ot a_1 \ot ...\ot a_n) = a_0 [D, a_1]...[D, a_n] \, . \label{eq15bis}
\end{equation}
One can then check the various conditions which in the commutative case
suffice to characterize Riemannian geometry \cite{Co3}, \cite{[Cncg]}.

\begin{theorem}\label{th3}
\begin{itemize}
\item[{\rm a)}] The spectral triple $(C^{\ify} (S^{4}_{\theta}),  \Hc , D)$
fulfills all axioms of noncommutative manifolds.
\item[{\rm b)}] Let $e \in C^{\ify} (S^{4}_{\theta} , M_4 (\Cb))$ be the
canonical idempotent given in (III.\ref{idempot}).
The Dirac operator $D$ fulfills
$$
\left\langle \left( e - \frac{1}{2} \right) [D,e]^{4} \right\rangle = \g
$$
where $\langle \ \rangle$ is the projection on the commutant of
$M_4 (\Cb)$ and $\g$ is the grading operator.
\end{itemize}
\end{theorem}

\bigskip
\noindent  The real structure \cite{[Coo2]} is given by the charge
conjugation operator $J$,
which involves in the noncommutative case the Tomita-Takesaki antilinear
involution.
The order one condition,
\begin{equation}
[[D,a],b^0] = 0 \qquad \fl \, a,b \in C^{\ify} (S^{4}_{\theta}) \, .
\label{eq:(4.102)}
\end{equation}
where $b^0 = J b^* J^{-1}$ follows easily from the derivation rules for the
$\delta_j$. \\
As we shall mention in the next section, Poincar\'e duality continues to hold.
\vspace{1cm}

\section{Isospectral deformations}

\noindent
We shall show in this section how to extend Theorem~\ref{th3} of the
previous section
to arbitrary metrics on the sphere $S^4$ which are invariant under rotation
of $u$ and $v$ and
have the same volume form as the one of the round metric.

\noindent
We shall in fact describe a very general construction of isospectral
deformations
of noncommutative geometries
which implies in particular that any
  compact spin Riemannian manifold $M$ whose isometry group has rank
$\geq 2$ admits a
natural one-parameter isospectral deformation to noncommutative geometries
$M_\theta$. The deformation of the algebra will be performed
  along the lines of \cite{[Ri1]}.

\noindent
We let $(\Ac , \Hc , D)$ be the canonical spectral triple associated with a
compact Riemannian spin manifold $M$. We recall that $\Ac = C^\ify(M)$ is
the algebra of smooth
functions on $M$, $\Hc= L^2(M,S)$ is the Hilbert space of spinors and $D$
is the Dirac operator.
We let $J$ be the charge conjugation operator which is an antilinear
isometry of $\Hc$.

\noindent
Let us assume that the group ${\rm Isom}(M)$ of isometries of $M$ has rank
$r\geq2$.
Then, we have an inclusion
\begin{equation}
\Tb^2 \subset {\rm Isom}(M) \, ,
\end{equation}
with $\Tb = \Rb / 2 \pi \Zb$ the usual torus, and we let $U(s) , s \in
\Tb^2$, be
the corresponding unitary operators in $\Hc = L^2(M,S)$ so that by construction
\begin{equation}
U(s) \, D = D \, U(s) \, , \, \, \, U(s) \, J = J \, U(s) \, .
\end{equation}
Also,
\begin{equation}
U(s) \, a \, U(s)^{-1} = \alpha_s(a) \, , \, \, \, \fl \, a \in \Ac \, ,
\label{actfun}
\end{equation}
where $\alpha_s \in {\rm Aut}(\Ac)$ is the action by isometries on the
algebra of functions on
$M$.

\noindent
We let $p = (p_1, p_2)$ be the generator of the two-parameters group $U(s)$
so that
\begin{equation}
U(s) = \exp(i(s_1 p_1 + s_2 p_2)) \, .
\end{equation}
The operators $p_1$ and $p_2$ commute with $D$ but anticommute with $J$.
Both $p_1$ and $p_2$
have integral spectrum,
\begin{equation}
{\rm Spec}(p_j) \subset \Zb \, , \, \, j = 1, 2 \, .
\end{equation}

\noindent
One defines a bigrading of the algebra of bounded operators in $\Hc$ with the
operator $T$ declared to be of bidegree
$(n_1,n_2)$ when,
\begin{equation}
\alpha_s(T) = \exp(i(s_1 n_1 + s_2 n_2)) \, T \, , \, \, \, \fl \, s \in
\Tb^2 \, ,
\end{equation}
where $\alpha_s(T) = U(s) \, T \, U(s)^{-1}$ as in (\ref{actfun}). \\
Any operator $T$ of class $C^\ify$ relative to $\alpha_s$ (i. e. such that
the map $s \rightarrow \alpha_s(T) $ is of class $C^\ify$ for the
norm topology) can be uniquely
written as a doubly infinite
norm convergent sum of homogeneous elements,
\begin{equation}
T = \sum_{n_1,n_2} \, \widehat{T}_{n_1,n_2} \, ,
\end{equation}
with $\widehat{T}_{n_1,n_2}$ of bidegree $(n_1,n_2)$ and where the sequence
of norms $||
\widehat{T}_{n_1,n_2} ||$ is of
rapid decay in $(n_1,n_2)$.

\bigskip

\noindent
Let $\lambda = \exp(2 \pi i \theta)$. For any operator $T$ in $\Hc$ of
class $C^\ify$ we define
its left twist $l(T)$ by
\begin{equation}
l(T) = \sum_{n_1,n_2} \, \widehat{T}_{n_1,n_2} \, \lambda^{n_2 p_1} \, ,
\end{equation}
and its right twist $r(T)$ by
\begin{equation}
r(T) = \sum_{n_1,n_2} \, \widehat{T}_{n_1,n_2} \, \lambda^{n_1 p_2} \, ,
\end{equation}
Since $|\lambda | = 1$ and $p_1$, $p_2$ are self-adjoint, both series
converge in norm. \\
One has,
\begin{lemma}\label{lem1}
\begin{itemize}
\item[{\rm a)}] Let $x$ be a homogeneous operator of bidegree $(n_1,n_2)$
and $y$ be
a homogeneous operator of  bidegree $(n'_1,n'_2)$. Then,
\begin{equation}
l(x) \, r(y) \, - \,  r(y) \, l(x) = (x \, y \, - y \, x) \,
\lambda^{n'_1 n_2} \lambda^{n_2 p_1 + n'_1 p_2}
\end{equation}
In particular, $[l(x), r(y)] = 0$ if $[x, y] = 0$.
\item[{\rm b)}] Let $x$ and $y$ be homogeneous operators as before and
define
\begin{equation}
x * y = \lambda^{n'_1 n_2} \, x y \, ; \label{star}
\end{equation}
then $l(x) l(y) = l(x * y)$.
\end{itemize}
\end{lemma}

\noindent
To check a) and b) one simply uses the following commutation rule which
is fulfilled for any
homogeneous operator $T$ of bidegree $(m,n)$,
\begin{equation}
\lambda^{a p_1 + b p_2} \, T = \lambda^{a m + b n} \, T \, \lambda^{a p_1 +
b p_2} \, ,
\, \,\, \forall a, b \in \Zb \, .
\end{equation}
One has then
\begin{equation}
l(x) \, r(y) = x \, \lambda^{n_2 p_1} \, y \, \lambda^{n'_1 p_2} = x \, y \,
\lambda^{n'_1 n_2} \, \, \lambda^{n_2 p_1 + n'_1 p_2} \,
\end{equation}
and
\begin{equation}
r(y) \, l(x) = y \, \lambda^{n'_1 p_2} \, x \, \lambda^{n_2 p_1} = y \, x \,
\lambda^{n'_1 n_2} \,\, \lambda^{n_2 p_1 + n'_1 p_2} \,
\end{equation}
which gives a). One checks b) in a similar way.

\noindent
The product $*$ defined in (\ref{star}) extends by linearity
to an associative product on the linear space of smooth operators and could
be called a $*$-product.

\noindent
One could also define a deformed `right product'. If $x$ is homogeneous of
bidegree
$(n_1,n_2)$ and $y$ is homogeneous of bidegree $(n'_1,n'_2)$ the product is
defined by
\begin{equation}
x *_{r} y = \lambda^{n_1 n'_2} \, x y \, .
\end{equation}
Then, along the lines of the previous lemma one shows that $r(x) r(y) = r(x
*_{r} y)$.

\bigskip
\noindent
Next, we twist the antiunitary isometry $J$ by
\begin{equation}
\wt{J} = J \, \lambda^{- p_1 p_2} \, . \label{realtwist}
\end{equation}
One has $\wt{J} = \lambda^{p_1 p_2} \, J $ and hence
\begin{equation}
\wt{J}^2  = J ^2 \, .
\end{equation}

\begin{lemma}\label{lem2}
For $x$ homogeneous of bidegree $(n_1,n_2)$ one has that
\begin{equation}
\wt{J} \, l(x) \, \wt{J}^{-1} = r(J \, x \, J^{-1}) \, \lambda^{-n_1 n_2} \, .
\end{equation}
\end{lemma}
For the proof one needs to check that
\begin{equation}
\wt{J} \, l(x) = r(J \, x \, J^{-1})  \, \lambda^{-n_1 n_2} \, \wt{J}.
\end{equation}
  One
has
\begin{equation}
\lambda^{-p_1 p_2} \, x = x \, \lambda^{-(p_1 + n_1) (p_2 + n_2)} = x \,
\lambda^{-n_1 n_2} \, \,
\lambda^{-(p_1 n_2 + n_1 p_2)} \, \,  \lambda^{-p_1 p_2}.
\end{equation}

\noindent  Then
\begin{equation}
\wt{J} \,
l(x) = J \, \lambda^{-p_1 p_2}
\, x \, \lambda^{n_2 p_1} = J
\, x \, \lambda^{-n_1 n_2} \,
\lambda^{- n_1 p_2} \, \, \lambda^{-p_1 p_2},
\end{equation}

\noindent  while

\begin{equation}
r(J \, x \, J^{-1}) \, \wt{J} = J \, x \,
J^{-1} \, \lambda^{-n_1 p_2} \,
J \, \lambda^{-p_1 p_2} = J \, x \, \lambda^{-n_1 p_2}  \,\, \lambda^{-p_1
p_2}.
\end{equation}
  Thus one gets the required equality.

\bigskip
\noindent
We can now define a new spectral triple where both $\Hc$ and the operator
$D$ are unchanged while the
algebra $\Ac$ and the involution $J$ are modified to $l(\Ac)$ and $\wt{J}$
respectively. By
Lemma~{\ref{lem1}}~b) one checks that  $l(\Ac)$ is still an algebra.

\noindent
Since $D$ is of bidegree $(0,0)$ one has,
\begin{equation}
[D, \, \l(a) ] = l([D, \, a]) \label{bound}
\end{equation}
which is enough to check that $[D, x]$ is bounded for any $x \in \l(\Ac)$.

\noindent
For $x,y \in l(\Ac)$ one checks that
\begin{equation}
[x, \, y^0] = 0 \, , \, \, \, y^0 = \wt{J} \, y^* \, \wt{J}^{-1} \, .
\label{commut}
\end{equation}
Indeed, one can assume that $x$ and $y$ are homogeneous and use
Lemma~{\ref{lem2}} together with Lemma~{\ref{lem1}}~a).
Combining equation (\ref{commut}) with equation (\ref{bound}) one then
checks the order one condition
\begin{equation}
[ \, [D, \, x ] \, , y^0 ] = 0 \, , \, \, \, \fl \, x, y \in \l(\Ac) \,
.\label{first}
\end{equation}

\noindent
As a first corollary of the previous construction we thus get
\begin{theorem}\label{th4}
Let $M$ be a compact spin Riemannian manifold whose isometry group has rank
$\geq 2$. Then $M$ admits a
natural one-parameter isospectral deformation to noncommutative geometries
$M_\theta$.
\end{theorem}
The deformed spectral triple is given by $(l(\Ac) , \Hc , D)$ with $\Hc =
L^2(M,S)$ the Hilbert space of spinors, $D$ the Dirac operator and $l(\Ac)$ is
really the algebra of smooth functions on $M$ with product deformed to the
$*$-product defined in (\ref{star}). Moreover, the real structure is given by
the twisted involution $\wt{J}$ defined in (\ref{realtwist}). One checks using
the results of \cite{Rieffel} and \cite{Co3} that Poincar\'e duality 
continues to hold for the
deformed spectral triple. We showed in
\cite{Co3} that the Dirac operator for the Levi-Civita connection
  minimizes the
action functional $ \int\!\!\!\!\!\!-  \, D^{2-n} $ (where $n$ is the 
dimension of $M$)
  among operators of the form $D+T$ which $\epsilon$ commute with $J$ 
and have the same commutators
  as $D$  with
any $a \in \Ac$ (so that $T$ belongs to the commutant of $\Ac$). It 
is important to check that this continues to hold
  in the deformed case. This is easy to see since we can also assume 
invariance under the action
$U(s) \, T \, U(s)^{-1} = \alpha_s(T)$ so that the space of available 
perturbations $T$ is
smaller in the deformed case.

\bigskip
\noindent
The above construction also allows us to extend Theorem~\ref{th3} of the
previous section
to arbitrary metrics on the sphere $S^4$ which are invariant under rotation
of $u$ and $v$ and
have as volume form $\sqrt{g} d x$ the round one.

\bigskip
\noindent In \cite{N-S} Nekrasov and Schwarz  showed
that
Yang-Mills gauge theory on
noncommutative $\Rb^4$ gives a conceptual understanding of the nonzero
B-field desingularization of the moduli space of instantons obtained by
perturbing
the ADHM equations \cite{[At]}. In  \cite{Witten}, Seiberg and
Witten
exhibited the unexpected relation
between the
standard gauge theory and the noncommutative one. The above work raises the
specific question
for NC-spheres $S^{4}_{\theta}$ whether one can implement such a
Seiberg-Witten relation as an
isospectral one. It also suggests to extend the above isospectral
deformations (Theorem~\ref{th4}) to more general compatible Poisson structures
on a given spin Riemmannian manifold.

\vspace{1cm}

\section{Final remarks}

\noindent
We shall end this paper with several important remarks,

\bigskip
\noindent {\it  The odd case}

\noindent  First there are formulas for the {\it odd } Chern character in
cyclic homology,
similar to those of section II above. Given an invertible element $ u \in
GL_r(\Ac)$,
the component ${\rm ch}_{n+\frac{1}{2}} (u)$ of its Chern character is as
above an
element of
\begin{equation}
\Ac \ot \underbrace{\ov{\Ac} \ot \cdots \ot \ov{\Ac}}_{2n-1} \label{eq12.2}
\end{equation}
where $\ov{\Ac} = \Ac / \Cb 1$ is the quotient of $\Ac$ by the scalar
multiples of
the unit 1.

\noindent The formula for ${\rm ch}_{n+\frac{1}{2}} (u)$ is (with $\lb_n$ a
normalization constant),
\begin{eqnarray}
{\rm ch}_{n+\frac{1}{2}} (u) &=& \lb_n  \{ \sum  u_{i_0 i_1}  \,
\ot u^{-1}_{i_1 i_2} \ot u_{i_2 i_3} \cdots \ot u^{-1}_{i_{2n-1} i_0} \\
&-& \sum  u^{-1}_{i_0 i_1}  \,
\ot u_{i_1 i_2} \ot u^{-1}_{i_2 i_3} \cdots \ot u_{i_{2n-1} i_0}
\} \,  \nonumber
\end{eqnarray}
As in the even case, the crucial property of the
components ${\rm ch}_{n+\frac{1}{2}} (u)$ is that they define a {\it cycle}
in the $(b,B)$
bicomplex of cyclic homology,
\begin{equation}
B \, {\rm ch}_{n-\frac{1}{2}} (u) = b \, {\rm ch}_{n+\frac{1}{2}} (u) \, .
\label{eq12.4}
\end{equation}
For any pair of integers $m,r$ we can define the odd analogues $\Bc_{m,r}$
as generated by the $r^2$ elements $u_{ij}$;
$i,j \in \{ 1, \ldots , r \}$ and we impose as above the relations
\begin{equation}
u \, u^*=u^* \, u = 1 \qquad u = [u_{ij}] \label{eq12.6}
\end{equation}
and
\begin{equation}
{\rm ch}_{j+\frac{1}{2}} (\rho(u)) = 0 \qquad \fl \, j < m \, . \label{eq2.7bis}
\end{equation}
One can prove as an exercice that the suspension of the corresponding NC
spaces
are contained in the
$Gr_{m,2r}$.

\bigskip

\noindent {\it  The Dirac operator and quantum groups}

\noindent There exists formulas for $q$-analogues of the Dirac
operator on
quantum groups,
(cf. \cite{biku}, \cite {[Majid2]}); let us call $Q$ these ``naive" Dirac
operators.
Now the fundamental equation to define the sought for true Dirac operator
$D$ which we used above implicitly on the deformed 3-sphere (after
suspension to the 4-sphere and for deformation parameters which are complex
of modulus one)
is,
\begin{equation}
  [D]_{q^2} = Q \, . \label{eq13.7}
\end{equation}
where the symbol $[x]_q$ has the usual meaning in $q$-analogues,
\begin{equation}
  [x]_q = \frac{q^x - q^{-x}}{q - q^{-1}} \, . \label{eq13.8}
\end{equation}
  The main point is that it is only
by virtue of this equation that the commutators $ [D, a] $ will be bounded,
and they will be so not only for the natural action of the algebra $\Ac$ of
functions
on $SU(2)_q$ on the Hilbert space of spinors but also for the natural action of
the opposite algebra $\Ac^o$; this is easy to prove in Fourier. But it is not
true
that $ [Q, a]$  is bounded, for $a \in \Ac $, due to the unbounded nature
of the
bimodule defining the $q$-analogue of the differential calculus.

\bigskip

\noindent {\it  Yang-Mills theory}

\noindent  One can develop the Yang-Mills theory on $S^{4}_{\theta}$ since
we now have all the
required structure, namely the algebra, the calculus and the ``vector
bundle" $e$ (naturally
endowed, in addition, with a prefered connection $\nabla$).
One can check that the basic results
of \cite{[Co]} apply. In particular Theorem 4, p 561 of \cite{[Co]}  gives
a basic inequality
showing that the Yang-Mills
  action, $YM(  \nabla) ={\int \!\!\!\!\!\! -} \,  \t ^2 \,  ds^{4} $,
(where $\t = \nabla^2$ is the curvature, and
$ds=D^{-1}$)
  has a strictly positive
  lower bound given by the topological invariant $\int\!\!\!\!\!- \g (e
-\frac{1}{2}) [D, e]^4
\ ds^4 = 1$.
\noindent  The next step is thus to extend the results of \cite{[At]} on the
classification of Yang-Mills connections to this situation. This was done
in \cite{Co-R}
  for the noncommutative torus and in \cite{N-S} for noncommutative $ \Rb^4$.
  Note however that in the noncommutative case the NC-sphere
  $S^{4}_{\theta}$ {\it is not isomorphic} to the
one-point compactification of noncommutative $ \Rb^4$ used there.
In particular, and in contrast to what happens for noncommutative $ 
\Rb^4$, even the
measure theory of
$S^{4}_{\theta}$ is very sensitive to the irrationality of the parameter
$\theta$.

\vfill\eject

\end{document}

\bibitem{[Ri1]} {\sc M.A. Rieffel}, {\it Deformation quantization for 
actions of R^d},
  Memoirs AMS 506 (1993).

\bibitem{[Majid]}  {\sc S. Majid}, {\it Foundations of Quantum Group Theory},
Cambridge University Press (1995).